\documentclass[journal,twoside,web]{ieeecolor}
\usepackage{setspace}
\usepackage[utf8]{inputenc}
\usepackage{color,soul}
\usepackage[T1]{fontenc}
\usepackage{multicol}
\usepackage{multirow}
\usepackage{textcomp}
\usepackage{mathrsfs}
\usepackage{graphicx}
\usepackage{epstopdf}
\usepackage{amssymb}
\usepackage{amsmath}
\usepackage{epstopdf}
\usepackage{psfrag}

\usepackage{amsthm}
\usepackage[tight,footnotesize]{subfigure}
\usepackage{cite}
\usepackage{generic}
\usepackage{enumerate}
\usepackage{mathtools}
\usepackage{mathabx}
\usepackage{algorithmic,algorithm}

\theoremstyle{plain}
\newtheorem{assumption}{Assumption}
\newtheorem{theorem}{Theorem}
\theoremstyle{definition}

\newtheorem{definition}{Definition}
\newtheorem{lemma}{Lemma}

\newtheorem{remark}{Remark}
\newtheorem{example}{Example}

\makeatletter
\newcommand{\vast}{\bBigg@{3.2}}
\newcommand{\Vast}{\bBigg@{5.5}}

\makeatother

\def\BibTeX{{\rm B\kern-.05em{\sc i\kern-.025em b}\kern-.08em
		T\kern-.1667em\lower.7ex\hbox{E}\kern-.125emX}}
\graphicspath{{Figure/}}
 
\begin{document}

\title{Optimal Decentralized Control for Uncertain\\ Systems by Symmetric Gauss-Seidel Semi-Proximal ALM}


\author{
	Jun Ma, 
	Zilong Cheng,
	Xiaoxue Zhang, 
	Masayoshi Tomizuka, \IEEEmembership{Life Fellow,~IEEE,}
	and Tong Heng Lee	
	\thanks{J. Ma is with the Department of Electrical and Computer Engineering, National University of Singapore, Singapore 117583, he is also with the Department of Mechanical Engineering, University of California, Berkeley, CA 94720 USA (e-mail: elemj@nus.edu.sg; jun.ma@berkeley.edu).}
	\thanks{Z. Cheng, X. Zhang, and T. H. Lee are with the NUS Graduate School for Integrative Sciences and Engineering, National University of Singapore, Singapore 119077 (e-mail: zilongcheng@u.nus.edu; xiaoxuezhang@u.nus.edu; eleleeth@nus.edu.sg).}
		\thanks{M. Tomizuka is with the Department of Mechanical Engineering, University of California, Berkeley, CA 94720 USA (e-mail: tomizuka@berkeley.edu).}
}

\maketitle

\begin{abstract}
	The
	$\mathcal{H}_2$ guaranteed cost decentralized control problem
	is investigated in this work.
	More specifically, on the basis of an appropriate $\mathcal{H}_2$ re-formulation that we put in place,
	the optimal control problem in the presence of parameter uncertainties is then suitably characterized by convex restriction and solved in parameter space.
	It is shown that a set of stabilizing decentralized controller gains
	for the uncertain system is parameterized in a convex set through appropriate convex restriction, 
	and then an approximated conic optimization problem is constructed. This facilitates the use of the symmetric Gauss-Seidel (sGS) semi-proximal augmented Lagrangian method (ALM),
	which attains high computational effectiveness. A comprehensive analysis is given on the application of the approach
	in solving the optimal decentralized control problem;
	and subsequently, the preserved decentralized structure, robust stability,
	and robust performance
	are all suitably guaranteed with the proposed methodology.
	Furthermore, an illustrative example is presented to demonstrate the effectiveness of the proposed optimization approach.
\end{abstract}


\section{Introduction}
Decentralized control has been widely applied to various large-scale systems
in view of several important attractive properties,
such as parallel computation of control variables,
ease of communication among controllers, sensors, and actuators, etc.
Typical applications of decentralized control in these situations
include electrical power network models,
transportation systems,
communication networks, robotics, etc.
Further, in general, a large-scale multi-input-multi-output (MIMO) system
can essentially be decomposed into a set of physically coupled subsystems,
and it then allows the independent design of controllers suitably.
As the decentralized control approach has a multi-level control architecture,
the control system thus also exhibits some good robustness properties
towards various structural perturbations that can occur in the high-level control function,
such as the breakdown of information exchanges between subsystems.

Stabilization of a MIMO system by decentralized control has been
rather thoroughly studied in the literature,
with methodologies
such as Lyapunov-type methods~\cite{bellman1962vector}
and algebraic-type methods~\cite{wang1973stabilization,davison1976robust}.
On the other hand,
although there are already some optimal and robust control techniques
that are relatively straightforward and well-established, it is still not easy to apply these relevant optimal and robust control techniques
to decentralized control;
such as the common-place usage of algebraic Riccati equations (AREs)
for linear quadratic regulator (LQR),
linear quadratic Gaussian (LQG), $\mathcal{H}_2$, and $\mathcal{H}_\infty$ control in the non-decentralized control situation.
The key evident reason is that the decentralized controller gain matrix exhibits
specific sparsity constraints.
In view of the sparsity constraints due
to this decentralized control architecture,
various alternate controller design methods
have been presented to cater to this shortcoming~\cite{vandenberghe2015chordal,fazelnia2016convex,wang2018separable,furieri2020learning,fattahi2019efficient}.
Thus in~\cite{witsenhausen1971separation},
through a problem re-formulation,
sufficient conditions are given
such that the standard LQG theory could be applied.
Elsewhere in~\cite{geromel1982optimal},
in the LQR and $\mathcal{H}_2$ control problems,
it is noted that
the gradient of a pre-defined objective function with respect
to the controller gain admits a closed-form solution,
which enables
the use of the gradient descent method (after the projection
of the gradient onto the decentralized constraint hyperplanes).
However, this approach only gives the local optimum
that is highly dependent on controller initialization,
and it does not take the robustness issue into consideration.

To further address
these shortcomings, several works on
parameterization have been proposed
where approaches are developed such that
the constrained optimization problem can be solved in extended parameter space.
For example in~\cite{geromel1991convex,geromel1994decentralized},
it is shown that
the sparsity constraints can be explicitly expressed in the parameter space,
thereby invoking the development and utilization of
a cutting-plane optimization algorithm.
In this framework,
the global optimum is then determined through an iterative framework,
and also robustness towards parameter uncertainties
is ensured~\cite{ma2019robust,ma2019parameter}.
Nevertheless, a very real practical limitation is that
the convergence of the cutting-plane optimization algorithm is rather slow.
The first reason is that,
due to the requirement in the procedure for the generation of separating hyperplanes,
additional equality constraints are added to the linear programming problem during each iteration
(thereby adding to the computational burden).
The second reason is that the nonlinear constraints are solved through outer linearization,
which is not computationally effective;
and these constraints are not always exhaustively exactly satisfied
but violated in a small number of instances, even upon completion of all iterations. Recent works on decentralized control by the use of block-diagonal Lyapunov functions are presented in~\cite{zheng2019distributed}, which ensure a property of sparsity invariance, and then the structural constraint on the controller can be suitably transformed into the constraints on Lyapunov variables. A complete generalization of this line of work is given in~\cite{furieri2019separable}, which utilizes a convex restriction to reformulate the unstructured problem to an equivalent convex program involving additional convex constraints to guarantee the structural constraint on the controller.

Added to all the above observations, quadratic invariance (QI) serves as the necessary and sufficient condition for the set of Youla parameters that respect the original information constraint to be convex~\cite{lessard2011quadratic, lessard2015convexity}, and then a feasible controller admitting the optimal performance can be determined. Similar research results and insights are also reported in~\cite{qi2004structured,rotkowitz2005characterization,furieri2019input,lin2019convex}. However, QI is a somewhat stringent condition for the convexity of decentralized control design problems. To overcome the limitation of QI, some further researches have been conducted. For instance, \cite{rotkowitz2011nearest} presents the closest subset and superset of the decentralized constraint, which are quadratically invariant when the original problem is not. The notion of sparsity invariance proposed in~\cite{furieri2020sparsity} can also be applied to design optimally distributed controllers subject to sparsity constraints on the controller structure; and it evidently outperforms the nearest QI subset in some cases and is also naturally applicable to the distributed static controller. However, these techniques are also computationally expensive, and the required conditions based on QI can also be rather too stringent and rigorous for many practical cases.

In more recent works,
a new optimization technique called the augmented Lagrangian method (ALM) has been presented,
which has attracted considerable attention from researchers
in the optimal control area~\cite{lin2013design,lee2017complexity}.
However, for some ill-conditioned optimization problems,
a slow rate of convergence can result from the application of the conventional ALM,
and for some large-scale optimization problems,
the conventional ALM cannot even be used to solve the given problem successfully.
Since the scale of the $\mathcal H_2$ optimization problem
under parameter uncertainties is usually very large,
and the condition number of the resulting optimization problem is typically not clear,
the conventional ALM is not feasible
as a proposed methodology for the specified optimization problem here.
To cater to the shortcomings of the conventional ALM,
in~\cite{li2016schur},
a semi-proximal term is introduced into the augmented Lagrangian function
to ensure that the related sub-problems can be solved effectively.
Additionally,
in~\cite{li2019block},
an approach denoted as
the symmetric Gauss-Seidel (sGS) technique
is applied to ensure the large-scale problem can be separated into a group of sub-problems.
Nevertheless, while promising,
all these methodologies essentially provide only rather generic guidelines at this stage,
and a more definite strategy and formulation (with accompanying comprehensive analysis)
is still lacking
on how to extend their usages in real-world problems.

With all of the above descriptions as a back-drop,
in this work, we develop and propose
a definite strategy and formulation that
uses the semi-proximal ALM and the sGS
to specifically solve the
optimal decentralized control problem with parameter uncertainties.
On the basis of a parameter space formulation with convex restriction,
the optimization problem
is restricted and approximated as a convex optimization problem. Then, the sGS semi-proximal ALM is utilized to determine the optimal solution effectively.
Importantly here in our formulation,
the sparsity constraints resulting from the decentralized control structure are satisfied.
Additionally, robust stability in the presence of parameter uncertainties is guaranteed,
and the system performance is maintained within a prescribed level.
This approach thus addresses all the shortcomings as previously highlighted,
and also attains high computational effectiveness
in solving
the optimal decentralized control problem.

The remainder of this paper is organized as follows.
In Section II, the problem formulation of the $\mathcal{H}_2$ guaranteed cost decentralized control problem
is given in its usual commonly-encountered form.
In Section III, we present and show how
the resulting optimization problem
can be suitably
constructed in the parameter space through convex restriction.
Next, since this problem now takes the structure of a convex optimization problem,
we then present and develop the specific procedures
utilizing the sGS semi-proximal ALM for effective optimization in this specific
parameter space
approximation of the $\mathcal{H}_2$ guaranteed cost decentralized control problem.
In Section V, an illustrative example is given to validate the results.
Finally, pertinent conclusions are drawn in Section VI.

\section{Problem Statement}
The following notations are used in the remaining text. $\mathbb R^{m\times n}$ ($\mathbb R^{n}$) denotes the real matrix with $m$ rows and $n$ columns ($n$ dimensional real column vector). $\mathbb S^{n}$ denotes the $n$ dimensional real symmetric matrix. $\mathbb S^{n}_{+}$ ($\mathbb S^{n}_{++}$) denotes the $n$ dimensional positive semi-definite (positive definite) real symmetric matrix. The symbol $A \succ 0$ ($A \succeq 0$) means that the matrix $A$ is positive (semi-)definite, and  $A\succ B$ ($A\succeq B$) means $A-B$ is positive (semi-)definite. $A^T$ ($x^T$) denotes the transpose of the matrix $A$ (vector $x$). $I_n$ represents the identity matrix with a dimension of $n\times n$. The operator $\operatorname{Tr}(A)$ denotes the trace of the square matrix $A$. The operator $\langle A, B \rangle$ denotes the Frobenius inner product i.e. $\langle A,B\rangle= \operatorname{Tr}\left(A^TB\right)$ for all $A,B \in \mathbb R^{m\times n}$. The norm operator based on the inner product operator is defined by $\|x\|=\sqrt{\langle x,x\rangle}$ for all $x\in \mathbb R^{m\times n}$. $\| H(s)\|_2$ represents the $\mathcal{H}_2$-norm of $H(s)$. $\otimes$ denotes the Kronecker product. $\operatorname{eig(A)}$ represents all the eigenvalues of the matrix $A$. $\operatorname{blocdiag}\{A_1, A_2, \cdots, A_n\}$ denotes a block diagonal matrix with diagonal entries $A_1, A_2, \cdots, A_n$. $\operatorname{diag}\{a_1, a_2, \cdots, a_n\}$ denotes a diagonal matrix with diagonal entries $a_1, a_2, \cdots, a_n$.

Consider a linear time-invariant (LTI) system
\begin{IEEEeqnarray}{rl}
	\dot x(t) &= A x(t) + B_2 u(t) + B_1 w(t), \label{eq:ssm1} \\
	z(t) &= C  x(t)+ D  u(t), \label{eq:ssm2}
\end{IEEEeqnarray}
with a static state feedback controller
\begin{IEEEeqnarray}{rl}
	u(t) = -K x(t), \label{eq:ssm4}
\end{IEEEeqnarray}
where $x\in \mathbb{R}^{n}$ is the state vector with $x(0)=x_0$, $u\in \mathbb{R}^{m}$ is the control input, $w\in \mathbb{R}^{l}$ is the exogenous disturbance input, $z\in \mathbb{R}^{q}$ is the controlled output,   $K\in \mathbb{R}^{m \times n}$ is the feedback gain matrix. $A$, $B_1$, $B_2$, $C$, and $D$ are constant real matrices with appropriate dimensions. It is assumed that there is no cross weighting between the state variables and the control variables, i.e. $C^T D=0$, and the control weighting matrix is nonsingular, i.e. $D^T D \succ 0$. Also, as a usual practice, it is assumed that $(A, B_2)$ is stabilizable and the pair $(A, C)$ has no unobservable modes on the imaginary axis. Remarkably, for a decentralized control problem, $K$ is constrained to be block diagonal.

Here, the objective function is defined as
\begin{IEEEeqnarray}{rl}
	J = \int_0^{\infty} z(t)^T z(t) \, dt.  \label{eq:Cost1}
\end{IEEEeqnarray}
To optimize \eqref{eq:Cost1}, it is equivalent to minimize the $\mathcal{H}_2$-norm of the transfer function
\begin{IEEEeqnarray}{rl}
	H(s)=(C- DK)(sI_n-A+B_2K)^{-1} B_1, \label{eq:Transfer function}
\end{IEEEeqnarray}
from $w$ to $z$, and the objective function \eqref{eq:Cost1} can be re-formulated as
\begin{IEEEeqnarray}{rl}
	J(K)&=\| H(s)\|_2^2 
	=\operatorname{Tr}\left((C-D K)W_c(C-DK)^T\right),  \label{eq:2_Cost2}
\end{IEEEeqnarray}
where $W_c$ is the controllability Gramian associated with the closed-loop system.

If all the associated matrices in the system are precisely known (there is no parameter uncertainty), the problem of finding the optimal decentralized $K$ to minimize the objective function \eqref{eq:Cost1} is considered as an $\mathcal{H}_2$ decentralized control problem. On the other hand, with the existence of parameter uncertainties, the design of a decentralized controller such that the upper bound to the $\mathcal{H}_2$-norm is minimized is referred to as an $\mathcal{H}_2$ guaranteed cost decentralized control problem.

\section{Optimization Problem Formulation in Parameter Space}
It is assumed that  $A$ and $B_2$  are subjected to parameter uncertainties. Define $p=n+m$, and then the following extended matrices are introduced for an alternative representation of the system:
\begin{gather}
	F=\begin{bmatrix}
		A & -B_2\\0  & 0
	\end{bmatrix} \in \mathbb{R}^{p\times p}, \quad G=\begin{bmatrix}
		0 \\I_m
	\end{bmatrix} \in \mathbb{R}^{p\times m},  \nonumber\\ Q=\begin{bmatrix}
		B_1 B_1^T &0 \\0 &0
	\end{bmatrix} \in \mathbb{S}^p,  \quad
	R=\begin{bmatrix}
		C^T C &0 \\0 & D^T D
	\end{bmatrix} \in \mathbb{S}^p.
\end{gather}

\begin{assumption}
	The parameter uncertainties are structural and convex-bounded.
\end{assumption}

Followed by Assumption 1, it is assumed that $F$ belongs to a polyhedral domain, which is expressed by a convex combination of the extreme matrices, where  $F=\sum_{i=1}^M \xi_i F_i$, $\xi_i \geq 0$, $\sum_{i=1}^M \xi_i=1$ , $F_i=\begin{bmatrix}
A_i & -B_{2,i}\\
0 & 0
\end{bmatrix} \in \mathbb{R}^{p \times p}$ denotes the extreme vertex of the uncertain domain. Remarkably, the precisely known system is a special case of the above expression, where $M=1$.

Denote the block partition of a matrix $W$ as
\begin{IEEEeqnarray}{rl}
	W=\begin{bmatrix}
		W_1 & W_2 \\
		W_2^T & W_3
	\end{bmatrix} \in \mathbb{S}^p,\label{eqn:partition} \end{IEEEeqnarray} with $W_1 \in \mathbb{S}^n_{++}$, $W_2 \in \mathbb{R}^{n \times m}$, $W_3 \in \mathbb{S}^{m}$, and then denote $ \Theta_i(W)= F_i W+W F_i^T+ Q$ as a block partitioned matrix, where
\begin{IEEEeqnarray}{rl}
	\Theta_i(W)= \begin{bmatrix}  \Theta_{i,1}(W) &  \Theta_{i,2}(W) \\ \Theta_{i,2}^T(W) &  \Theta_{i,3}(W) \end{bmatrix}\in \mathbb{S}^p,
\end{IEEEeqnarray}
with $\Theta_{i,1}(W) \in \mathbb{S}^{n}, \Theta_{i,2}(W) \in \mathbb{R}^{n \times m}, \Theta_{i,3}(W) \in \mathbb{S}^{m}$. With an appropriate convex restriction, Theorem~\ref{thm:1} presents a subset of controller gains that preserve the decentralized structure, robust stability, and robust performance, in the presence of parameter uncertainties.

\begin{theorem}\label{thm:1}
	For the controller with a specific decentralized structure
	\begin{IEEEeqnarray}{rl}~\label{eq:decen}
		K  = \operatorname{blocdiag}\left\{K_\textup{D,1}, K_\textup{D,2}, \cdots, K_{\textup{D},m}\right\},
	\end{IEEEeqnarray}
	with $K_{\textup{D},j} \in \mathbb{R}^{{1\times \textup{D}_j}}$, $\forall j=1,2, \cdots, m$, one can define the set 
	\begin{IEEEeqnarray}{lll}
		\mathscr{C}  = \Big\{&W\in\mathbb S^p: W\succeq0,    \Theta_{i,1}(W)  \preceq 0,\forall i=1,2,\cdots,M,\IEEEnonumber\\
		&W_\textup{1}=\operatorname{blocdiag}\{W_\textup{1D,1}, W_\textup{1D,2},\cdots, W_{\textup{1D},m}\},\IEEEnonumber\\ &W_\textup{2}=\operatorname{blocdiag}\{W_\textup{2D,1}, W_\textup{2D,2}, \cdots, W_{\textup{2D},m}\},\IEEEnonumber\\
		&W_\textup{1D,j} \in \mathbb{R}^{{\textup{D}_j \times \textup{D}_j}},W_\textup{2D,j} \in \mathbb{R}^{\textup{D}_j},\forall j =1,2,\cdots,m\Big\},\IEEEeqnarraynumspace 
	\end{IEEEeqnarray}
	and 
	\begin{IEEEeqnarray}{lll}
		\mathscr{K}= \left\{K=W_2^T W_1^{-1}: W \in \mathscr{C}\right\}, \label{eq:add}
	\end{IEEEeqnarray}
	 then,
	\begin{enumerate}[(a)]
		\item $K \in \mathscr{K}$ holds the decentralized structure in \eqref{eq:decen}.
		\item $K \in \mathscr{K}$ stabilizes the closed-loop system in the presence of model uncertainties.
		\item  $K \in \mathscr{K}$ gives $\langle R,W\rangle \geq \| H_i(s)\|_2^2, \forall i=1,2,\cdots,M$ where $\|H_i(s)\|_2$ represents the $\mathcal H_2$-norm with respect to the $i$th extreme system.
	\end{enumerate}
\end{theorem}

\noindent \textbf{Proof of Theorem \ref{thm:1}:}
For Statement (a), it can be easily derived that
$
W_2^T=\operatorname{blocdiag}\{W_\textup{2D,1}^T,W_\textup{2D,2}^T,\cdots, W_\textup{2D,m}^T\}
$
and
$
W_1^{-1} = \operatorname{blocdiag}\{W_\textup{1D,1}^{-1},W_\textup{1D,2}^{-1},\cdots,W_\textup{1D,m}^{-1}\}
$.
Then from~\eqref{eq:add}, 
$K=  \operatorname{blocdiag}\Big\{W_\textup{2D,1}^T W_\textup{1D,1}^{-1}, W_\textup{2D,2}^T W_\textup{1D,2}^{-1},  \cdots, W_\textup{2D,m}^T W_\textup{1D,m}^{-1}\Big\},$
which holds the decentralized structure in~\eqref{eq:decen}.

For Statement (b), $\Theta_{i,1}(W) \preceq 0$ is equivalent to
\begin{IEEEeqnarray}{rl}
	A_i W_1-  B_{2,i} W_2^T+W_1  A_i^T-W_2   B_{2,i}^T+B_1 B_1^T   \preceq 0. 
\end{IEEEeqnarray}
Since $W_1\succ 0$, we have
\begin{IEEEeqnarray}{rCl}~\label{eq:c}
	(  A_i-  B_{2,i} W_2^T W_1^{-1})W_1+W_1 (  A_i &-  B_{2,i} W_2^T W_1^{-1})^T  \nonumber\\
	&+B_1 B_1^T  \preceq 0.
\end{IEEEeqnarray}
It is straightforward that $K=W_2^T W_1^{-1}$ stabilizes the extreme system and then we have
$W_2=W_1  K^T$, and subsequently we can construct $W= \begin{bmatrix} W_1 & W_1 K^T\\KW_1 &W_3 \end{bmatrix}$. From Schur complement, $W_3$ is a free variable to choose such that  $W \succeq 0$ is ensured. $W \succeq 0$ is not a necessary condition to stabilize the extreme systems, but it indeed provides a norm bound to the controller gain. It is straightforward that, to ensure the stability over the entire uncertain domain which is convex, it suffices to check the stability at the vertices of the convex polyhedron. Therefore, if the stability holds for all the extreme systems, then the stability for the entire uncertain domain is guaranteed.

For Statement (c), first, we have $W_1 \succeq W_{c,i}$. Then, from Schur complement, $W\succeq 0$ leads to
\begin{IEEEeqnarray*}{rCl}
	W_3 &\succeq& W_2^T W_1^{-1}W_2
	= KW_1  K^T
	\succeq  KW_{c,i}  K^T.\IEEEyesnumber
\end{IEEEeqnarray*}
Therefore, it follows that
$
	\langle R,W\rangle
	\geq \| H_i(s)\|_2^2.\IEEEyesnumber
$ \hfill{\qed}

\begin{definition}
	The system~\eqref{eq:ssm2} is called robustly strongly decentralized stabilizable
	if $\mathscr{C}\neq \emptyset$.
\end{definition}

It can be seen that $\langle R,W\rangle$ provides an upper bound to $\| H_i(s)\|_2^2$. Then, it is aimed to solve the optimization problem $W = \operatorname{argmin} \{\langle R,W\rangle: W \in \mathscr{C}\}$, which yields $K=W_2^{T} W_1^{-1} \in \mathscr{K}$, such that the upper bound to the $\mathcal{H}_2$-norm is minimized. Hence,  the optimization problem is summarized in the following form:
\begin{IEEEeqnarray*}{rl}\label{eq:opt1}
	\displaystyle\operatorname*{minimize}_{W\in\mathbb S^{p}}\quad &  \langle R,W\rangle  \\
	\operatorname*{subject\ to}\quad &W \succeq 0, \\
	&\Theta_{i,1}(W)   \preceq 0, \, \forall i=1,2,\cdots,M \\
	&W\in\Phi(s),\IEEEyesnumber
\end{IEEEeqnarray*}
and $\Phi(s)$ denotes the set of all $W$ satisfying the sparsity constraints.

\begin{remark}
	Problem~\eqref{eq:opt1} is a convex restriction of the original $\mathcal{H}_2$ guaranteed cost decentralized control problem, and Theorem~\ref{thm:1}
	characterizes a subset of stabilizing decentralized controller gains. Pertinent approaches of convex restriction can be found in~\cite{geromel1991convex,geromel1994decentralized}. Besides, several approaches can be used to reduce the conservatism, such as those presented in~\cite{furieri2019separable,furieri2020sparsity}.
\end{remark}

It is obvious that the operator $\Theta_{i,1}(W)$ is a bounded linear operator (affine to $W$), and the sparsity constraints are linear equality constraints. Since the objective function is also a linear function, and the optimization variable $W$ is confined in a convex cone, thus it falls into the category of the convex optimization problem, where $\mathscr{C}$ is a convex set. In particular, the optimization problem is a linear conic programming problem. In this paper, the sGS semi-proximal ALM is introduced to solve the given optimization problem.

To express the optimization problem explicitly, we define a matrix
$
V=\begin{bmatrix}
I_n & 0_{n \times m}
\end{bmatrix}
$, and then the optimization problem can be equivalently expressed in the matrix form, where
\begin{IEEEeqnarray*}{rl}\label{equation:optimization_problem1}
	\displaystyle\operatorname*{minimize}_{W\in\mathbb S^{p}}\quad& \langle R,W\rangle\IEEEnonumber\\
	\operatorname*{subject\ to}\quad
	&W\in \mathbb S^{p}_+\\
	&-V\left( F_i W+W  F_i^T+Q\right)V^T \in \mathbb S^{n}_+\\
	&-V_{j1}WV_{j2}=0\\
	&\forall i=1,2,\cdots,M,\,\forall j=1,2,\cdots,N.\IEEEyesnumber
\end{IEEEeqnarray*}

\begin{remark}\label{property1}
	By splitting $W_1$ and $W_2$ into $m^2$ sub-blocks, the zero blocks in $W_1$ can be expressed by $m(m-1)/2$ equality constraints (because $W_1$ is symmetric), and the zero blocks in $W_2$ can be expressed by another $m(m-1)$ equality constraints. Therefore, $N=3m(m-1)/2$. Obviously, adjacent zero blocks can be combined into one single equality constraint, thus $N=3m(m-1)/2$ is not the minimum number of the equality constraints. For the sake of illustration purposes and without loss of generality, the adjacent zero blocks are not combined in the following analysis.
\end{remark}

Define $\Psi_i = -V( F_i W+W  F_i^T+Q)V^T,\forall i=1,2,\cdots,M$. Then the optimization problem can be denoted in a compact form, which is shown as
\begin{IEEEeqnarray*}{rl}\label{equation:optimization_problem2}
	\displaystyle\operatorname*{minimize}_{W\in\mathbb S^{p}}\quad
	& \langle R,W\rangle\\
	\operatorname*{subject\ to}\quad
	&
	\mathcal G(W) \in\mathcal K,  \IEEEyesnumber
\end{IEEEeqnarray*}
where $\mathcal G(W)$ is a linear mapping given by
\begin{IEEEeqnarray}{l}
	\mathcal G(W)=
	(W,\Psi_1,\Psi_2,\cdots,\Psi_M,\IEEEnonumber\\
	\qquad\qquad V_{11}WV_{12},V_{21}WV_{22},\cdots,V_{N1}WV_{N2}),
\end{IEEEeqnarray}
and the convex cone $\mathcal K$ can be denoted as
\begin{IEEEeqnarray*}{rCl}
	\mathcal K&=&\mathbb S^{p}_+\times\underbrace{\mathbb S^{n}_+\times \cdots\times \mathbb S^{n}_+}_M\\
	&&\times\underbrace{\{0_{v_{11}\times v_{12}}\}\times\{0_{v_{21}\times v_{22}}\}\times \cdots\times \{0_{v_{N1}\times v_{N2}}\}}_N, \IEEEeqnarraynumspace\yesnumber
\end{IEEEeqnarray*}
where for all $i=1,2,\cdots,N$, the scalars $v_{i1}$ and $v_{i2}$ denote the number of rows of the matrix $V_{i1}$ and the number of the columns of the matrix $V_{i2}$, respectively. Since the positive semi-definite cone is self-dual, it is straightforward to express the dual cone $\mathcal K^*$ as
\begin{IEEEeqnarray*}{rCl}
	\mathcal K^*&=&\mathbb S^{p}_+\times\underbrace{\mathbb S^{n}_+\times \cdots\times \mathbb S^{n}_+}_M\\
	&&\times\underbrace{\mathbb R^{v_{11}\times v_{12}}\times\mathbb R^{v_{21}\times v_{22}}\times \cdots\times \mathbb R^{v_{N1}\times v_{N2}}}_N.\IEEEyesnumber
\end{IEEEeqnarray*}
Define the linear space $\mathcal X$ in terms of the cone $\mathcal K$, which is given by
\begin{IEEEeqnarray*}{rCl}
	\mathcal X&=&\mathbb S^{p}\times\underbrace{\mathbb S^{n}\times \cdots\times \mathbb S^{n}}_M\\
	&&\times\underbrace{\mathbb R^{v_{11}\times v_{12}}\times\mathbb R^{v_{21}\times v_{22}}\times \cdots\times \mathbb R^{v_{N1}\times v_{N2}}}_N,\IEEEyesnumber
\end{IEEEeqnarray*}
and it is straigtforward that the cone $\mathcal K\subset \mathcal X$ and the dual cone $\mathcal K^*\subset \mathcal X$.


\begin{assumption}\label{assumption:strong_duality}
	Problem~\eqref{eq:opt1} is strictly feasible.
\end{assumption}

Under Assumption~\ref{assumption:strong_duality},  Slater's condition is satisfied. Therefore,  strong duality always holds for the proposed optimization problem, and the optimal solution to the linear conic optimization problem~\eqref{equation:optimization_problem2} can be obtained by solving the corresponding dual problem, due to the difficulty to deal with the primal problem directly. In terms of the optimization problem~\eqref{equation:optimization_problem2}, the Lagrangian function is defined as
\begin{IEEEeqnarray*}{rCl}
	\mathcal L(W;X)&=& \langle R,W\rangle-\left\langle X,\mathcal G(W)\right\rangle,\IEEEyesnumber
\end{IEEEeqnarray*}
where $X=(X_0,X_1,\cdots,X_{M+N}) \in\mathcal K^*$ is the Lagrange multiplier.
It follows that the Lagrangian dual function $\theta(X)$ is obtained by
\begin{IEEEeqnarray*}{rCl}
	\theta(X)
	&=& \displaystyle\min_{W\in\mathbb S^{p}} \Big\{\langle R,W\rangle-\langle X,\mathcal G(W)\rangle\Big\}.\IEEEyesnumber
\end{IEEEeqnarray*}
It is shown that the Lagrangian dual function can be denoted in the explicit form, where
\begin{IEEEeqnarray*}{l}
	\displaystyle\min_{W\in\mathbb S^{p}} \Big\{\langle R,W\rangle-\langle X,\mathcal G(W)\rangle\Big\}\\
	\begin{array}{l}
		=\left\langle X_1+X_2+\cdots+X_M,VQV^T\right\rangle+\\
		\displaystyle\min_{W\in\mathbb S^{p}}
		\left\{\begin{array}{c}
			\langle R,W\rangle - \langle X_0,W\rangle\\
			+\displaystyle\sum_{i=1}^M \left\langle  F_i^TV^TX_iV+V^TX_iVF_i  ,W\right\rangle\\
			+\frac{1}{2}\displaystyle\sum_{j=1}^N \left\langle V_{j1}^TX_{M+j}V_{j2}^T+V_{j2}X_{M+j}^TV_{j1},W\right\rangle
		\end{array}\right\}.
	\end{array}\\
	\IEEEyesnumber
\end{IEEEeqnarray*}
For the sake of simplicity in the remaining text, we define
\begin{IEEEeqnarray*}{l}
	\mathcal F(X_1,X_2,\cdots,X_M)=-\left\langle X_1+X_2+\cdots+X_M,VQV^T\right\rangle\\
	\mathcal A(X_0,X_1,\cdots,X_{M+N})=-X_0+ \sum_{i=1}^M\Big(F_i^TV^TX_iV\\
	+V^TX_iV F_i\Big)+\frac{1}{2}\sum_{j=1}^N\Big(V_{j1}^TX_{M+j}V_{j2}^T+V_{j2}X_{M+j}^TV_{j1}\Big).
	\IEEEyesnumber
\end{IEEEeqnarray*}
For all $i=0,1,\cdots,M+N$, $\mathcal A_i (X_0, X_1, \cdots,X_{i-1},X_{i+1},$ $X_{i+2},\cdots, X_{M+N})$ is used to denote the linear functions, in which the terms related to $X_i$ is removed from the linear function $\mathcal A(X_0,X_1\cdots,X_{M+N})$. 


Now we can denote the Lagrangian dual function explicitly, which is given by
\begin{IEEEeqnarray*}{l}
	\theta(X)=
	\begin{cases}
		-\mathcal F(X_1,X_2,\cdots,X_M)&\text{if }\mathcal A(X)+R=0\\
		-\infty &\text{otherwise}.
	\end{cases}\IEEEeqnarraynumspace\IEEEyesnumber\\
\end{IEEEeqnarray*}
The Lagrangian dual problem is to maximize the Lagrangian dual function, and thus, the Lagrangian dual problem is given by
\begin{IEEEeqnarray*}{rl}
	\operatorname*{minimize}_{X \in\mathcal X}\quad& \mathcal F(X_1,X_2,\cdots,X_M)\\
	\operatorname*{subject\ to}\quad & \mathcal A(X)+R=0,\, X\in\mathcal K^*.\IEEEyesnumber
\end{IEEEeqnarray*}

For any linear space $\mathcal Y$ and any convex set $\mathcal C\subset \mathcal Y$, define the indicator function $\delta_\mathcal C(v)$ such that for any $v\in\mathcal Y$, it follows that
\begin{IEEEeqnarray*}{rl}
	\delta_\mathcal C(v)=\begin{cases}
		0 & \text{if } v\in\mathcal C\\
		+\infty  & \text{otherwise.}
	\end{cases}\IEEEyesnumber
\end{IEEEeqnarray*}
Finally, the optimization problem can be equivalently converted to the form with only one linear equality constraint, which is given in the following form:
\begin{IEEEeqnarray*}{rl}\label{equation:dual_problem_1}
	\operatorname*{minimize}_{X \in\mathcal X}\quad& \mathcal F(X_1,X_2,\cdots,X_M)+\delta_{\mathbb S_+^p}(X_0)+\delta_{\mathbb S_+^n}(X_1)\\
	&+\delta_{\mathbb S_+^n}(X_2)+\cdots+\delta_{\mathbb S_+^n}(X_M)\\
	\operatorname*{subject\ to}\quad & \mathcal A(X)+R=0.\IEEEyesnumber
\end{IEEEeqnarray*}

\section{sGS Semi-Proximal ALM for Optimal Decentralized Control}
In the following text, the sGS semi-proximal ALM is presented to solve the dual problem~\eqref{equation:dual_problem_1}. Notably, in the remaining text in this section, variable $X$ is considered as the optimization variable of the problem~\eqref{equation:dual_problem_1} and the variable $W$ is considered as the Lagrange multiplier.

\noindent{\textbf{Step 1: Initialization.}}

Choose the parameters $\sigma>0$ and $\tau\in(0,(1+\sqrt 5)/2)$, the parameter $\alpha_i$ such that the linear operator $\mathcal S_i$ is a positive operator, the initial matrices $(X^0;W^0)\in\mathcal X\times \mathbb S^{p}$, and the parameter of the stopping criterion $\epsilon>0$.

\noindent{\textbf{Step 2: Update of the optimization variable $\textit{X}$.}}

{{For backward sGS sweep,}}
define the augmented Lagrangian function as
\begin{IEEEeqnarray*}{rCl}
	\mathcal L_\sigma(X;W)&=&
	\mathcal F(X_1,X_2,\cdots,X_M)+\delta_{\mathbb S_+^p}(X_0)+\delta_{\mathbb S_+^n}(X_1)\\
	&&+\delta_{\mathbb S_+^n}(X_2)+\cdots+\delta_{\mathbb S_+^n}(X_M)\\
	&&+\frac{\sigma}{2}\left\|\mathcal A(X)+R-\sigma^{-1}W\right\|^2-\frac{1}{2\sigma}\|W\|^2,\IEEEyesnumber
\end{IEEEeqnarray*}
where $W\in\mathbb S^{p}$ is the Lagrange multiplier for the augmented Lagrangian function. Notice that there are three groups of variables in the optimization problem: the variable $X_0$ is related to the positive semi-definite constraint for the Lagrange multiplier $W$; the variables $X_1,X_2,\cdots,X_M$ are regarding the positive semi-definite constraints of a group of linear functions; the variables $X_{M+1},X_{M+2}\cdots,X_{M+N}$ are related to the linear equality sparsity constraints for the Lagrange multiplier $W$. The sub-problems in the backward sGS sweep will be solved in terms of these three groups of variables separately.

We begin with the sub-problem of the third group. For all $i=N,N-1,\cdots,1$, the sub-problem with respect to the variable $X_{M+i}$ is given by
\begin{IEEEeqnarray*}{rCl}\label{equation:sub_problem_1}
	\bar X_{M+i}^{k+1}&=&\displaystyle\operatorname*{argmin}_{X_{M+i}\in\mathbb R^{v_{i1}\times v_{i2} }} \quad \mathcal L_\sigma\Big(X_0^{k},X_1^{k},\cdots,X_{M+i-1}^k,\\
	&&X_{M+i},\bar X_{M+i+1}^{k+1},\bar X_{M+i+2}^{k+1},\cdots ,\bar X_{M+N}^{k+1};W^k\Big).\IEEEeqnarraynumspace\IEEEyesnumber
\end{IEEEeqnarray*}
where $\bar X_{M+i}^{k+1}$ denotes the solution to the sub-problem~\eqref{equation:sub_problem_1} in the $(k+1)$th iteration, $X_0^{k},X_1^{k},\cdots,X_{M+i-1}^k$, and $W^k$ denote the values of the optimization variable and the Lagrange multiplier in the $k$th iteration, respectively, and $\bar X_{M+i+1}^{k+1},\bar X_{M+i+2}^{k+1},\cdots ,\bar X_{M+N}^{k+1}$ denote the updated optimization variables in the $(k+1)$th iteration.

Since the sub-problem is an unconstrained optimization problem, the optimality condition is given by
\begin{IEEEeqnarray*}{rCl}
	0
	&=&  V_{i1}\Big[\mathcal A\Big(X_0^{k},X_1^{k},\cdots,X_{M+i-1}^k,X_{M+i},\\
	&&\bar X_{M+i+1}^{k+1},\bar X_{M+i+2}^{k+1},\cdots ,\bar X_{M+N}^{k+1}\Big)+R-\sigma^{-1}W^k\Big]V_{i2},\\\IEEEyesnumber
\end{IEEEeqnarray*}
and then we have
\begin{IEEEeqnarray*}{rCl}\label{equation:sub_problem_2}
	\bar X_{M+i}^{k+1} &=&  -2V_{i1}\Big[\mathcal A_{M+i}\Big(X_0^{k},X_1^{k},\cdots,X_{M+i-1}^k,\bar X_{M+i+1}^{k+1},\\
	&&\bar X_{M+i+2}^{k+1},\cdots ,\bar X_{M+N}^{k+1}\Big)+R-\sigma^{-1}W^k\Big]V_{i2}.\IEEEyesnumber
\end{IEEEeqnarray*}


In terms of minimizing the augmented Lagrangian function with respect to the variable $X_i$ for all $i=M, M-1,\cdots,1$, a proximal term must be considered in the sub-problem during the iterations. To include the proximal term without influencing the convergence of the algorithm, we firstly introduce the positive linear operator. For any linear space $\mathcal Y$, a linear operator $\mathcal S:\mathcal Y\rightarrow \mathcal Y$ is positive, which means that for all $v\in\mathcal Y$, it follows that $\langle v,\mathcal S(v)\rangle\ge 0$.
The sub-problem in terms of the variable $X_i$ in the backward sGS sweep is given by
\begin{IEEEeqnarray*}{rCl}
	\bar X_i^{k+1}&=&\operatorname*{argmin}_{X_i\in\mathbb S^{n}} \quad \mathcal L_\sigma\Big(X_0^{k},X_1^{k},\cdots,X_{i-1}^k,X_{i},\bar X_{i+1}^{k+1},\\
	&&\bar X_{i+2}^{k+2},\cdots ,\bar X_{M+N}^{k+1};W^k\Big)+\frac{1}{2}\left\|X_i-X_i^k \right\|_{\mathcal S_i}^2,\IEEEeqnarraynumspace\IEEEyesnumber
\end{IEEEeqnarray*}
where $\mathcal S_i$ is a positive linear operator. To solve this sub-problem effectively,   $\mathcal S_i$ is chosen as
\begin{IEEEeqnarray*}{rCl}
	\mathcal S_i(X)&=&\alpha_i \mathcal I(X)-\sigma V F_iF_i^TV^TX-\sigma V F_iV^TXV F_iV^T\\
	&& -\sigma V F_i^TV^TXVF_i^TV^T-\sigma XV F_iF_i^TV^T,\IEEEyesnumber
\end{IEEEeqnarray*}
where $\mathcal I$ denotes an identity operator, and one positive choice of $\alpha_i$ is the maximum eigenvalue of the vectorization matrix, which is
\begin{IEEEeqnarray*}{l}~\label{eq:alphai}
	\alpha_i = \sigma \max \big\{\operatorname{eig}\big(VF_i F_i^TV^T\otimes I + I \otimes VF_iF_i^TV^T\\+ (VF_i^TV^T)\otimes (VF_iV^T)
	+(VF_iV^T)\otimes(VF_i^TV^T)\big)\big\}.\IEEEeqnarraynumspace\IEEEyesnumber
\end{IEEEeqnarray*}

Before presenting the optimality condition to the sub-problem,  Theorem~\ref{theorem:projection_operator} is introduced to determine the  projection operator.
\begin{theorem}\label{theorem:projection_operator}
	The projection operator $\Pi_\mathcal C(\cdot)$ with respect to the convex cone $\mathcal C$ can be expressed as
	\begin{IEEEeqnarray}{rCl}
		\Pi_{\mathcal C}=(I+\alpha \partial \delta_{\mathcal C})^{-1},
	\end{IEEEeqnarray}
	where $\partial (\cdot)$ denotes the sub-differential operator, and $\alpha \in \mathbb R$ can be an arbitrary real number.
\end{theorem}

\noindent \textbf{Proof of Theorem \ref{theorem:projection_operator}:}
Define a finite dimensional Euclidean space $\mathcal X$ equipped with an inner product and its induced norm such that $\mathcal C\subset \mathcal X$. For any $x\in\mathcal X$, there exists $z\in\mathcal X$ such that $z\in(I+\alpha \partial \delta_{\mathcal C})^{-1}(x)$. Then it follows that
\begin{IEEEeqnarray}{rCl}\label{equation:projection_lemma}
	x\in(I+\alpha \partial \delta_{\mathcal C})(z)=z+\alpha \partial \delta_\mathcal C(z).
\end{IEEEeqnarray}
Note that the projection operator $\Pi_\mathcal C(z)$ can be expressed as
\begin{IEEEeqnarray}{rCl}\label{equation_projection_lemma2}
	\Pi_\mathcal C(z)=\operatorname*{argmin}_{z\in W}\left\{\delta_{\mathcal C}(x)+\frac{1}{2\alpha}\|z-x\|^2\right\}.
\end{IEEEeqnarray}
Since the optimization problem in \eqref{equation_projection_lemma2} is strictly convex, the sufficient and necessary optimality condition for the optimization problem of the projection operator can be expressed as
\begin{IEEEeqnarray}{rCl}
	0\in\alpha\partial \delta_\mathcal C(x)+z-x,
\end{IEEEeqnarray}
which is equivalent to \eqref{equation:projection_lemma}. Note that the projection onto a convex cone is unique. Therefore, the mapping from $x$ to $z$ is also unique, which means the operator $(I+\alpha \partial \delta_\mathcal C)^{-1}(\cdot)$ is a point-to-point mapping.   \hfill{\qed}

Then the optimality condition to the sub-problem of the second group is given by
\begin{IEEEeqnarray*}{rCl}
	0
	&\in&-V QV^T
	+\partial \delta_{\mathbb S_+^n}(X_i)\\
	&&+\sigma V F_i\Big[\mathcal A\Big(X_0^{k},X_1^{k},\cdots,X_{i-1}^k,X_{i},\bar X_{i+1}^{k+1},\\
	&&\;\;\;\bar X_{i+2}^{k+1},\cdots ,\bar X_{M+N}^{k+1}\Big)+R-\sigma^{-1}W^k\Big]V^T\\
	&&+\sigma V\Big[\mathcal A\Big(X_0^{k},X_1^{k},\cdots,X_{i-1}^k,X_{i},\bar X_{i+1}^{k+1},\\
	&&\;\;\;\bar X_{i+2}^{k+2},\cdots ,\bar X_{M+N}^{k+1}\Big)+R-\sigma^{-1}W^k\Big] F_i^TV^T\\
	&&+\mathcal S_i\left(X_i-X_i^k\right).\IEEEyesnumber
\end{IEEEeqnarray*}
Define
\begin{IEEEeqnarray*}{l}
	\text{LHS}_i\Big(X_0^{k},X_1^{k}\cdots,X_{i}^k,\bar X_{i+1}^{k+1},\bar X_{i+2}^{k+1},\cdots ,\bar X_{M+N}^{k+1};W^k\Big)\\
	=VQV^T-\sigma V F_i\Big[\mathcal A_i\Big(X_0^{k},X_1^{k}\cdots,X_{i-1}^k,\bar X_{i+1}^{k+1},\\
	\quad\quad\;\bar X_{i+2}^{k+1},\cdots ,\bar X_{M+N}^{k+1}\Big)+R-\sigma^{-1}W^k\Big]V^T\\
	\quad-\sigma V\Big[\mathcal A_i\Big(X_0^{k},X_1^{k},\cdots,X_{i-1}^k,\bar X_{i+1}^{k+1},\\
	\quad\quad\;\bar X_{i+2}^{k+1},\cdots ,\bar X_{M+N}^{k+1}\Big)+R-\sigma^{-1}W^k\Big] F_i^TV^T\\
	\quad-\Big(\sigma V F_iF_i^TV^TX_i^k+\sigma V F_iV^TX_i^kV F_iV^T\\
	\quad\quad+\sigma V F_i^TV^TX_i^kVF_i^TV^T+\sigma X_i^kV F_iF_i^TV^T-\alpha_i X_i^k\Big),\\ \IEEEyesnumber
\end{IEEEeqnarray*}  
and then it follows that
\begin{IEEEeqnarray*}{rCl}\label{eq:xx}
	\bar X_i^{k+1}&=&\Pi_{\mathbb S_+^n}\Big[
	\alpha_i^{-1}\text{LHS}_i\Big(X_0^{k},X_1^{k},\cdots,X_{i}^k,\\
	&&\bar X_{i+1}^{k+1},\bar X_{i+1}^{k+2},\cdots ,\bar X_{M+N}^{k+1};W^k\Big) \Big],\IEEEyesnumber
\end{IEEEeqnarray*}
where $\Pi_{\mathbb S_+^{n}}(\cdot)$ denotes the projection operator in terms of the positive semi-definite cone, and the projection results can be obtained by Lemma~\ref{lemma:proj}.
\begin{lemma}\label{lemma:proj}
	Projection onto the positive semi-definite cone can be computed explicitly. Let $X = \sum_{i=1}^n \lambda_i v_i v_i^T\in\mathbb S^n$ be the eigenvalue decomposition of the matrix $X$ with the eigenvalues satisfying $\lambda_1 \geq \lambda_2 \geq \dots \geq \lambda_n$, where $v_i$ denotes the eigenvector corresponding to the $i$th eigenvalue. Then the projection onto the positive semi-definite cone of the matrix $X$ can be expressed by
	\begin{IEEEeqnarray}{rCl}
		\Pi_{\mathbb S^n_+} (X) = \sum_{i=1}^n \max\left\{\lambda_i, 0\right\} v_i v_i^T.
	\end{IEEEeqnarray}
\end{lemma}

Then it is aimed to determine the solution to the sub-problem of the first group, where the sub-problem in terms of the variable $X_0$ in the backward sGS sweep is given as
\begin{IEEEeqnarray*}{l}
	X_0^{k+1}=\operatorname*{argmin}_{X_0\in\mathbb S^{p}} \mathcal L_\sigma\Big(X_0,\bar X_1^{k+1},\bar X_1^{k+2},\cdots,\bar X_{M+N}^{k+1};W^k\Big).\\\IEEEyesnumber
\end{IEEEeqnarray*}
The optimality condition is given by
\begin{IEEEeqnarray*}{rCl}
	0
	&\in&\sigma^{-1}\partial \delta_{\mathbb S_+^p}(X_0)+X_0\\
	&&-\Big[\mathcal A_0\Big(\bar X_1^{k+1},\bar X_2^{k+1},\cdots,\bar X_{M+N}^{k+1}\Big)+R-\sigma^{-1}W^k\Big],\\\IEEEyesnumber
\end{IEEEeqnarray*}
which gives
\begin{IEEEeqnarray*}{rCl}\label{equation:sub_problem_3}
	X_0^{k+1}&=&\Pi_{\mathbb S_+^p}\Big[\mathcal A_0\Big(\bar X_1^{k+1},\bar X_2^{k+1},\cdots,\bar X_{M+N}^{k+1}\Big)\\
	&&+R-\sigma^{-1}W^k \Big].\IEEEyesnumber
\end{IEEEeqnarray*}

{{For forward sGS sweep,}} for all $i=1,2,\cdots,M$, the second group optimization variables in the forward sGS sweep are given by
\begin{IEEEeqnarray*}{rCl}\label{equation:sub_problem_4}
	X_i^{k+1}&=&\Pi_{\mathbb S_+^n}\Big[
	\alpha_i^{-1}\text{LHS}_i\Big(X_0^{k+1},X_1^{k+1},\cdots,X_{i-1}^{k+1},\\
	&&\bar X_{i}^{k+1},\bar X_{i+1}^{k+1},\cdots ,\bar X_{M+N}^{k+1};W^k\Big) \Big].\IEEEyesnumber
\end{IEEEeqnarray*}
For all $i=1,2,\cdots,N$, the third group optimization variables in the forward sGS sweep are given by
\begin{IEEEeqnarray*}{l}\label{equation:sub_problem_5}
	X_{M+i}^{k+1} =  -2V_{i1}\Big[\mathcal A_{M+i}\Big(X_0^{k+1},X_1^{k+1},\cdots,X_{M+i-1}^{k+1},\\
	\bar X_{M+i+1}^{k+1},\bar X_{M+i+2}^{k+1},\cdots ,\bar X_{M+N}^{k+1}\Big)+R-\sigma^{-1}W^k\Big]V_{i2}.\IEEEeqnarraynumspace\IEEEyesnumber
\end{IEEEeqnarray*}

\noindent{\textbf{Step 3: Update of the Lagrange multiplier $\textit{W}$.}}

The Lagrange multiplier can be determined by
\begin{IEEEeqnarray}{rCl}\label{equation:ADMM_W}
	W^{k+1}=W^k-\tau\sigma\Big[\mathcal A\Big(X_0^{k+1},X_1^{k+1},\cdots,X_{M+N}^{k+1}\Big)+R\Big]. \IEEEeqnarraynumspace
\end{IEEEeqnarray}

\noindent{\textbf{Step 4: Stopping criterion.}}

The relative residual error in terms of the optimization variable $X_0$ is given by
\begin{IEEEeqnarray}{rCl}\label{equation:residual_error1}
	\text{err}_{X_0}^{k}&=&\frac{\left\|X_0^{k}-\Pi_{\mathbb S_{+}^{p}}(X_0^{k}-W^{k})\right\|}{1+\left\|W^{k}\right\|+\left\|X_0^{k}\right\|}.
\end{IEEEeqnarray}
For all $i=1,2,\cdots,M$, the relative residual errors in terms of the optimization variable $X_i$ are given by
\begin{IEEEeqnarray}{rCl}\label{equation:residual_error2}
	\text{err}_{X_i}^{k}&=&\frac{\left\|
		\begin{array}{l}
			X_i^k-\Pi_{\mathbb S_+^n}\Big[V QV^T+V F_iWV^T\\
			\quad\quad\quad\quad+VW F_i^TV^T+X_i^k\Big]
		\end{array}\right\|}{1+\left\|
		\begin{array}{c}
			V QV^T+V F_iWV^T\\
			+VW F_i^TV^T
		\end{array}
		\right\|+\left\|X_i^{k}\right\|}.
\end{IEEEeqnarray}
For all $i=1,2,\cdots,N$, the relative residual errors in terms of the optimization variable $X_{M+i}$ are given by
\begin{IEEEeqnarray}{rCl}\label{equation:residual_error3}
	\text{err}_{X_{M+i}}^{k}&=&\left\|V_{i1}W^kV_{i2}\right\|, 
\end{IEEEeqnarray}
and the relative residual error in terms of the equality constraint is given by
\begin{IEEEeqnarray}{rCl}\label{equation:residual_error4}
	\text{err}_{eq}^k&=&\frac{\left\|\mathcal A\left(X^k_0,X^k_1,\cdots,X^k_{M+N}\right)+R\right\|}{1+\left\|R\right\|}.
\end{IEEEeqnarray}
Define \begin{IEEEeqnarray}{rCl}\label{equation:maxresidual_error}
	\text{err}_{\max}=\max\left\{\text{err}_{X_0}^{k+1},\cdots,\text{err}_{X_{M+N}}^{k+1},\text{err}_{eq}^{k+1}\right\},
\end{IEEEeqnarray}
the optimization process will be terminated if $\text{err}_{\max} <\epsilon$, and $W^{k+1}$ is the optimal solution to the optimization problem~\eqref{equation:optimization_problem2}.

\noindent{\textbf{Step 5: Increase of the iteration index.}}

Set $k=k+1$ and go back to Step 2.



\section{Illustrative Example}
To illustrate the effectiveness of the above results, Example~\ref{exam:1} is given, which presents a decentralized controller design problem, where the state matrix and the control input matrix are randomly chosen such that their elements are stochastic variables uniformly distributed over $[0, 1]$.  The optimization algorithm is implemented in Python 3.7.5 with Numpy 1.16.4, and executed on a computer with 16G RAM and a 2.2GHz i7-8750H processor (6 cores). The parameters for initialization is given by: $\sigma=1$, $\tau=0.618$, $X^0=0$, $W^0=0$, and the stopping criterion is set as $\epsilon = 10^{-3}$.

\begin{example}\label{exam:1}
	Consider $x=\begin{bmatrix}
	x_1 & x_2 & x_3
	\end{bmatrix}^T$ and a linear system
	\begin{IEEEeqnarray*}{rl}
		\dot {{x}}&=A {x} +B_2 u+B_1 w,\\
		z&=Cx+Du,\\
		u&=-Kx,\nonumber
	\end{IEEEeqnarray*}
	where
		\begin{gather*}\nonumber
	A =
	\begin{bmatrix}
	0.1054&	0.6248&	0.1958\\
	0.2393&	0.6948&	0.6950\\
	0.4520&	0.3189&	0.8708
	\end{bmatrix}, \, B_1 =\begin{bmatrix}
	1 & 0 & 0\\
	0 & 1& 0\\
	0 &0 & 1
	\end{bmatrix},\\
	B_2  =
	\begin{bmatrix}
	0.9315	&0.7939\\
	0.9722	&0.1061\\
	0.5317	&0.7750		
	\end{bmatrix}, \,     C=\begin{bmatrix}
	1 & 0 & 0  \\
	0 &0 &0\\
	0  &0 &0
	\end{bmatrix}, \,
	D=\begin{bmatrix}
	0 & 0 \\1 & 0 \\0&1
	\end{bmatrix}, \nonumber
	\end{gather*}
	and $K$ is a block diagonal matrix with a prescribed structure.
\end{example}


To validate the effectiveness of the proposed methodology in terms of different scales, we vary the number of uncertain parameters in $A$ from 1, 4, to 9 by imposing parameter uncertainties to the leading principal sub-matrices of $A$ and also $A$ itself, with a magnitude of $\pm5\%$ of their nominal values. In this case, the number of extreme systems $M= 2$, 16, and 512, respectively. For $B_2$, no parameter uncertainty is imposed. 

Based on this example, the performance of our proposed methodology is compared with the solvers SCS and CVXOPT, and the results of comparison are shown in Table~\ref{table:tab}, where ``$\textup{Yes}$'' or ``$\textup{No}$'' denote the success or failure upon execution of the approach, respectively. It is observed that SCS is not able to return the solution successfully in all three cases. For CVXOPT, the solution can be obtained for the cases when $M=2$ or 16. However, it cannot scale to the case when $M=512$ and return the optimization results successfully. On the contrary, our proposed approach is capable of returning reliable optimization results successfully for all the scenarios. 
Also, when comparing the second-order algorithm to the first-order algorithm, it would be the case that when a relatively large duality gap is given, the first-order method should be faster theoretically, but the second-order method can achieve a very precise solution (which would be usually impractical for the first-order optimization algorithm). This is a generally typical observation that may be noted for such methodologies.

Taking the case that $M=512$ as an example, the optimal $W$ is given by
\begin{IEEEeqnarray*}{l}
	W=	\begin{bmatrix}
		9.5658  & -0.7779 & 0    &	1.2594 & 0\\
		-0.7779 & 0.5043 & 0     & 0.7632 & 0\\
		     0  & 0      & 1.5605 &0      & 4.6878  \\
	 	 1.2594 & 0.7632 & 0     &1.8644   & -0.0007 \\
	       	0 & 0        & 4.6878  & -0.0007 & 14.3733
	\end{bmatrix},
\end{IEEEeqnarray*}
and the optimal decentralized controller gain $K$ is given by
\begin{IEEEeqnarray*}{l}
	K=	 \begin{bmatrix}
	0.2913 & 1.9626 & 0\\
	   	0  & 0 & 3.0040
	\end{bmatrix}.
\end{IEEEeqnarray*}
\begin{table}[t!]\centering
	\footnotesize
	\caption{Comparison results of SCS, CVXOPT, and ALM} 
	\label{table:tab}
	\begin{tabular}{|c|c|c|c|}
		\hline
		& SCS      & CVXOPT   & ALM     \\ \hline
		M=2   & $\textup{No} $ & $\textup{Yes}$  & $\textup{Yes}$ \\ \hline
		M=16  & $\textup{No} $ & $\textup{Yes}$  & $\textup{Yes}$ \\ \hline
		M=512 & $\textup{No}  $ & $ \textup{No} $ & $\textup{Yes}$ \\ \hline
	\end{tabular}
\end{table}

\begin{figure}
	\centering
	\includegraphics[trim=10 50 10 80,width=0.80\columnwidth]{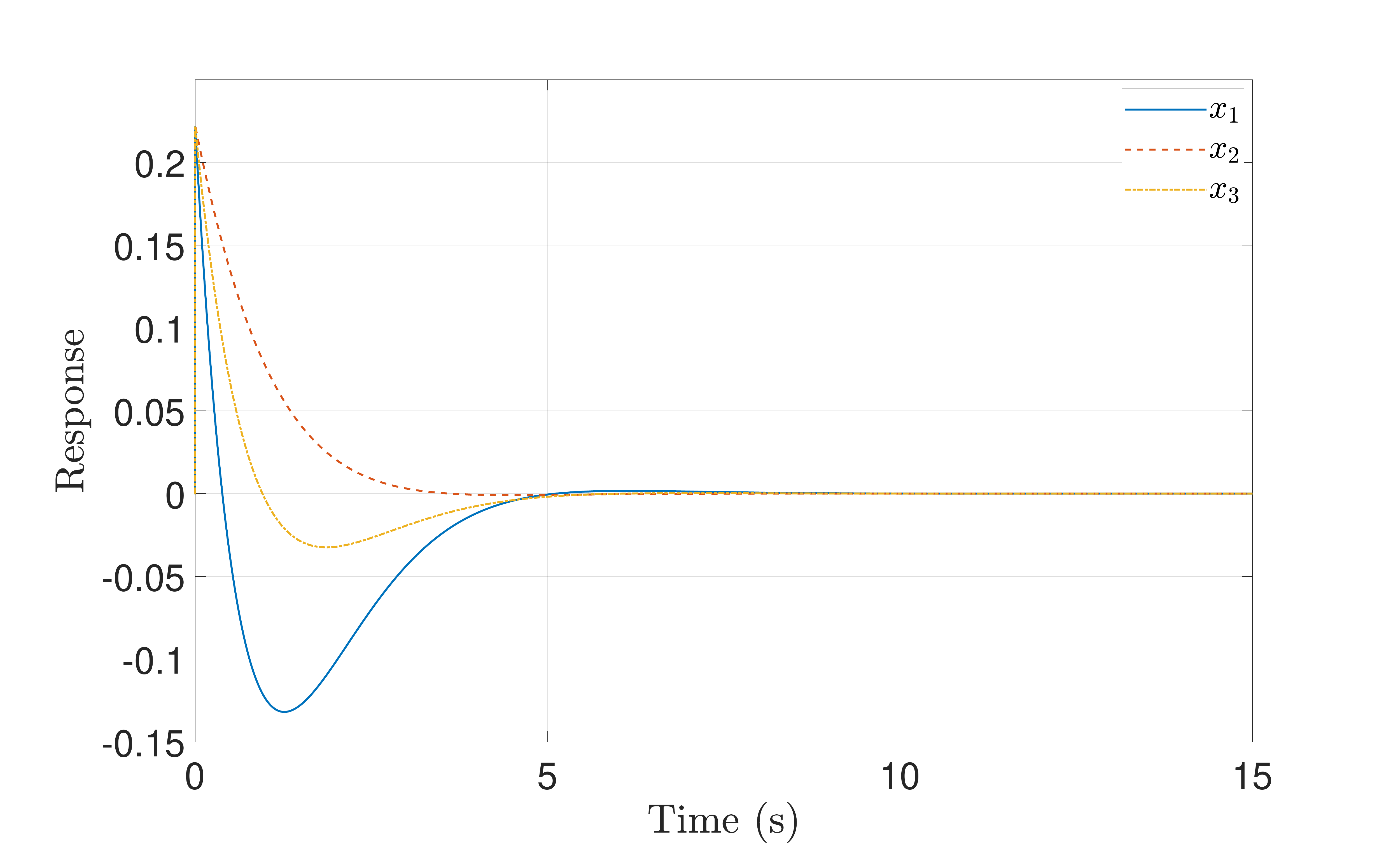}
	\caption{System response in Example 1}
	\label{fig:response1}
\end{figure}

It can be seen that the decentralized structure is preserved, and the upper bound to $\|H(s)\|_2^2$ is $\operatorname{Tr}(RW)=11.4302$. In the simulation, $w$ is characterized as a vector of the impulse disturbance, and an extreme system with all the uncertain parameters at their lower bounds is used. The responses of all  state variables are illustrated in Fig.~\ref{fig:response1}, and it can be seen that the robust stability is guaranteed for this extreme system. Also, it is worthwhile to mention that the stability is ensured for all the uncertain systems within the uncertain domain.   

\section{Conclusion}
This paper presents a highly effective
optimization algorithm for the optimal decentralized control problem with parameter uncertainties.
With parameterization of the stabilizing controller gain in the extended parameter space,
the problem is restricted and approximated as a convex optimization problem through appropriate convex restriction,
which is solved by the sGS semi-proximal ALM.
Furthermore, in the detailed development and methodology that we present here (which goes beyond just ``generic'' guidelines),
both the closed-loop stability
and the system performance are guaranteed in the presence of model uncertainties. An illustrative example demonstrates the applicability and effectiveness of the proposed approach.

\bibliographystyle{IEEEtran}
\bibliography{IEEEabrv,Reference}

\end{document}